# Naïve Philosophical Foundations (1967)

Bhupinder Singh Anand[1]

This soliloquy outlines some naïve philosophical arguments underlying the thesis that mathematics ought to be viewed simply as a universal set of languages, some of precise expression, and some of unambiguous and effective communication.

**Contents**

1. A Program

2. Guidelines

3. Do we know the source of the contradictions?

4. On what there is: In ontology and in language

5. An illustrative model: language and ontology

6. Logicism: Is the Russell-Frege definition of number significant?

7. Summary

## 1. A Program (1967)

One way of explaining the seeming incompatibilities amongst the strictures / structures provided by the intuitionistic, logicist, formalist, conventionalist and nominalist schools may be to regard these as guidelines not for the entire discipline accepted by common intuition as mathematics, but merely for varying particular aspects of such a concept.

One is then led to develop and isolate from these part-models the latter characteristics, and to in corporate them in a theory which may reflect a more comprehensive sub-structure of mathematics.

Such a view owes its conception in part to the development in physics of the quantum theory, which accommodated both wave and particle features in a mathematical - if not observably physical - form, whereas the earlier classical wave and particle models for light were contradictory.

---
[1] The author is an independent scholar. E-mail: re@alixcomsi.com; anandb@vsnl.com. Postal address: 32, Agarwal House, D Road, Churchgate, Mumbai - 400 020, INDIA. Tel: +91 (22) 2281 3353. Fax: +91 (22) 2209 5091.

A synthesised view of mathematics today should, it seems, be able to reflect in some form the basic issues on which the various schools were founded. Such as the logicist's identity of mathematics and logic, the formalist's stress on the internal validity and self-sufficiency criteria of a theory, the intuitionist's objection to passing from the negation of a general statement to an existential one without additional safeguards, the conventionalist's contention that the rules of a language delineate its ontology, as also the nominalist's scruples about the existence of classes of classes, amongst others.

The relative strength of such a construction would lie of course in the amount of mathematics containable in it, as also in the extent to which it can explicitly show how the various schools are illustrations of - models for - its sub-theories.

But would this not necessitate an extensive re-examination of current ideas and trends in philosophy and psychology as well as in logic?

## 2. Guidelines (1967)

(*a*) Are the contradictions arrived at fairly?

(*b*) Is the Russell-Frege definition of number significant?

(*c*) Why should arbitrarily formed predicates have a number?

(*d*) Are definitions true statements in the language or in the meta-language?

(*e*) Is the insistence on finitary means justified?

(*f*) Is the logicist thesis tenable?

(*g*) What is the position of the Law of the Excluded Middle?

(*h*) Does it make sense to say that logic has a domain?

(*i*) Are not the axioms of predicate logic tautologies?

(*j*) Doesn't predicate logic presuppose an understanding of sets?

(*k*) How can we give an open formula a truth-value that is not a set?

(*l*) What then is the nature of the connectives between open formulas?

(*m*) Do they merely clarify relations between existing sets?

(*n*) How can they be used to create new formulas for set construction?

(*o*) If we know '*f*(*x*) & *g*(*x*)' as true, and if *b* is an individual and *f*(*b*) is true, is it at all significant to say that by substitution, '*f*(*b*) & *g*(*b*)' being true, *g*(*b*) is deduced true?

(*p*) Doesn't the truth of *g*(*b*) follow immediately?

(*q*) What I am enquiring is whether the concept of deduction is at all meaningful, except as a mnemonic aid?

(*r*) What then is logic?

(*s*) Just as propositional logic presumes the content of its propositions - which are facts - so also may not predicate logic presume the content of its sentences - sets of facts?

(*t*) Doesn't logic merely direct sentences into possible categories?

(*u*) Are statements containing logical connectives on a different level from thos not containing them?

(*v*) What are many-valued logics on this view?

### 3. Do we 'know' the source of the contradictions? (20/03/1967)

(Cf. Quine, W.V. "On what there is". In "From a logical point of view". Harvard University Press. 1953.)

(*a*) Consider the expression:

    (*i*) $x \notin x$.

If we suppose that there is a class '*a*' whose members are precisely those that satisfy (*a*)(*i*), then we would hold that, in this instance, we have discovered a true statement schema:

    (*ii*) $x \in a$ if, and only if, $x \notin x$,

which expresses a host of facts concerning '*a*' and all the various members of some pre-existing universe.

But this belief is surely mistaken, for:

    (*iii*) $a \in a$ if, and only if, $a \notin a$,

is clearly false.

(*b*) Suppose, on the other hand, we say that we are defining a class '*a*' by:

   (*i*) $x \in a$ if, and only if, $x \notin x$.

Though this should be a true statement in our language about the theory, it may no longer be a statement in the language of the theory. But if we treat definition as a creative activity for producing a larger ontology, it is not surprising that we can arrive back at a supposedly true statement:

   (*ii*) $a \in a$ if, and only if, $a \notin a$,

inside the language. This position regarding creativity may differ but formally from our earlier platonistic stand.

(*c*) However, if we do not view definition as mere name-giving to newly born or already flourishing objects, then it is not easy to see what all the fuss is about.

For, if definition requires eliminability, then expressions such as '$a \in a$' and '$a \notin a$' are immediately suspect - since we are able to eliminate only '$x \in a$' from any expression.

And '*a*' in isolation is merely a strange creature giving rise to pseudo-expressions which confuse us as to their allowability into our language because of their familiar appearance.

But then, so too does Pegasus confuse us into sometimes creating worlds of ideas and unactualised possibilities!

And Quine has forcibly presented the case that a name need not name anything. For names belong to language essentially. And are easy to construct.

(*d*) There is a fuss, for the contradictions still haunt us some. So possibly we are loth to admit an error in our earliest 'discovery'. The seemingly true statement schema:

   (*i*) $x \in a$ if, and only if, $x \notin x$.

Now could it be that this reluctance to accept the negation of Cantor's Comprehension Axiom is psychologically motivated?

The cause to which we are clinging so stubbornly - armed with Russell's types, Zermelo's efforts, amongst others - may be that starting from an ontological stand of precise individuals and properties, we must somehow have the right to build up further properties into our universe. The paradoxes seem to prevent us from doing so with complete freedom.

(*e*) But why do we not feel the need to a similar liberty in the other direction? Regarding individuals.

Why do we not feel as strongly or as readily that by defining all the properties that occur in our ontology for a new individual, we may enlarge our universe?

The path may not be any smoother. For suppose we intend to introduce the individual '*k*' into our ontology. And our ontology contains a property schema '$P(x, y)$'. (Which may, for example be '$y$ loves $x$').

If our desire for liberty was sincere, we should feel free to then assign properties at will to the new entry.

But what happens?

Let us assign the $P(x, y)$'s to the entity '*k*' as follows:

(*i*)  $P(x, k)$ if, and only if, $\sim P(x, x)$.

Since '*k*' is part of our ontology, do we have:

(*ii*) $P(k, k)$ or $\sim P(k, k)$?

(*f*) My point is that as long as you have the desire to construct new relations amongst existing entities, you should also have the equal desire to construct new entities out of existing relations.

That if you have the feeling you can discover all kinds of possible relations amongst the individuals, you should also feel you can discover all kinds of individuals enmeshed in your relations.

That the guidelines in one case should be as useful in the other. That if every open formula in individuals seems to define a predicate, then every open formula in predicates should define an individual. To take a very naive view.

That we may be psychologically misled into feeling that a predicate open formula defines an entity known as the predicate of a predicate.

(*g*) So maybe there is much to be said for the nominalist stand. And isn't the idea that every individual be equivalnt to the set of all the predicates that it satisfies at the heart of Leibniz's notion of indiscernibles? As also at the heart of phenomenalism and positivism?

And where the external world is concerned, is it possible that quantum-interpreted phenomena may contain instances of plurality where the objects are indiscernibles - notwithstanding Leibniz's contention?

And inspite of Russell's claim of having no content to his universe does not the fact that it has no indiscernibles give it content - at least in the form of a special characteristic? Or is this too ambitious a claim?

## 4. On what there 'is': In ontology and in language (23/03/1967)

(*a*) I form concepts. That much seems reasonably clear to me. Their location I assume to be in the commonly referred to intuition. Concept space may be a better name for it.

(*b*) An analysis of these concepts I find to be a more difficult task than indicating their significance. So I intend to study merely the latter. However, I do take individuals, properties and facts as concepts.

(*c*) Events in physical space, indeed the space itself, are perceived and digested by my senses, whence they transform into concepts. Positivists would possibly claim that all my concepts are so derived.

(*d*) My concepts I may map into a language. This map you may decode into your concepts.

Assuming that both of us accept a common external world, I can understand why language is so useful.

(*e*) When I set up a language, there is what I talk about. Serious dispute cannot arise so long as my language faithfully refers to my concepts.

(*f*) I may feel the need to include Pegasus among my concepts. Your stoutest efforts will not convince me to analyse the name out a la Russell. A description into non-trivial terms of my ontology I would consider inadequate. And the trivial description of 'pegasises' I would only agree to as an introduction of a name for a concept of being Pegasus - a concept antecedent to the being of Pegasus among my concepts.

Or I may protest altogether against the being of any 'pegasises' concept in my concept space, and platonically refuse to admit discovery or creation of any such concept.

(*g*) Confusion may sometimes arise. You may wrongly translate my language into your concepts. My conceptual scheme may contradict the external world. I may have concepts not accessible to you.

In the first case you would be mistaken. In the second I should be convicted of error - and possibly idealism! But who is to judge?

Of some interest is the third. This I see as the cause of all genuine ontological disputes. From philosophy through to theology.

Taken to be a question of individual concepts, ontology seems more a matter of taste, inclination and, above all, feeling and belief in this case.

So its interest as a problem is, after all, trivial. As it should be.

(*h*) For, as long as I concern myself with ontology, restricting myself to a language constructed on the basis of my concepts, I shall for all practical purposes be dealing with the small aspect of the world which is conceptualised by my senses. And this, as Zeno's reflections seem to indicate, hardly can be said to exhaust nature's complexity.

(*i*) And if mathematics is to serve us in working with the real world, no satisfactory, or rather complete, mathematical conceptual scheme can be constructed only on intuitive concepts of the natural integers.

In fact, despite intuitionistic efforts, no adequate conceptual framework seems constructible to me.

(*j*) So I turn my back for the moment on concepts. All I am left with then is language, and possibly codifications of nature into language.

And my inability to grasp the totality of nature's concepts is contained in my use of variable names, and the transition from propositions to schemata.

And the test of any codifications as suitable for nature will be the inclusion in it of the concepts that are within my grasp.

(*k*) But what there 'is' in addition may, after all, depend on language in cases where empirical verification is lacking.

(*l*) And the language should presuppose some content - or at least possible content - and the logicist thesis could be misconceived. Or merely misrepresented.

## 5. An illustrative model: language and ontology (24/03/1967)

I have a concept of a possible universe that I should like to codify into language.

In my universe there are individuals, and there are properties. The landscape is otherwise deserted.

The individuals I shall name $a$, $b$, $c$, $d$, and $e$. The properties $F$, $G$, and $H$.

There are also (in some sense of being which is not entirely clear to me) facts in my universe. These I shall represent in my language as:

$F(a)$, $F(b)$, $G(b)$, $G(c)$, $G(e)$, $H(b)$, $H(c)$ and $H(e)$.

I shall call these true expressions in my language.

There are no such things (or whatever it is that facts are supposed to be) as non-facts in my universe. All the same, I admit certain expressions into my language - possibly for the sake of symmetry, but more so because tradition seems to demand such an action. These are:

$F(c)$, $F(d)$, $F(e)$, $G(a)$, $G(d)$, $H(a)$, and $H(d)$.

I shall call these false expressions.

Though my language, containing these expressions, is thus two-valued, in my universe there are only facts.

A very natural question may be asked for any set of individuals. Is there a property satisfied by all the members of the set, and none others?

I think I must be very clear about the nature of my enquiry. I am not asking whether my language can countenance the introduction of a further expression purporting to be a property. Such an entry, like the introduction of false expressions, may not present formidable difficulties. But I am enquiring whether my universe already contains such a property.

Taking $\{a, b, d\}$, as the set, I find no property which gives rise to true expressions for this set only. My finding is, of course, empirical.

For the set $\{a, b\}$ however, the property $F$ does give rise to true expressions; and no other individual satisfies $F$. And I may conveniently identify the set with F insofar as they are both names of the same entity.

What of the set $\{b, c, e\}$? Both $G$ and $H$ express facts for the members of this set only. But there is no unique property identifiable with this set. And, in passing, I may remark that such an event does not cause any concern usually. Properties with the same extension are tolerated easily.

I conclude that not every set of individuals can be identified with a unique property.

So, a set of individuals may not name anything in my universe.

A question of far greater significance is as to the nature of sets of properties. Classically these have been treated as being identifiable with a different quality of being in the universe from that of properties and individuals.

But though my language is prolific in sets, my universe is starved for entities. So I look for some more direct identifications for these sets than those suggested by precedent.

Surprisingly, I am successful - or so it seems. And my solution appears so natural - at least from nominalistic standards - that I begin to suspect that tradition may well have been merely disguising it.

For a set of properties, I ask the question whether any individual has just those properties, and none others.

For the set $\{F, G\}$ there is no such individual.

The set $\{F, G, H\}$ may be identified with the individual $b$, which is the only one satisfying all three properties. (I note, incidentally, that such identification has positivistic overtones.)

Similarly, $\{F\}$ may be identified with $a$.

But now I consider the set $\{G, H\}$. Both $c$ and $e$ satisfy only this set. Which is a most surprising characteristic of my universe. It contains two indiscernibles!

(Inspite of Leibniz, and Russell's subsequent backing of his ideas on the intuitive notion of equality, modern physics has made a universe with such characteristics rather feasible. What is required for such a feature is that some set of properties be identified with a plurality of individuals.)

I find, then, that not every set of properties is identifiable with an individual.

So, if I contain myself to the ontology outlined, some sets of properties, as also of individuals, don't exist, while some do, and still others exhibit an ambiguous character.

But all this is peculiar to my universe. And not every universe need be of this type.

The universe being constructed by an intuitionist may have differing qualities. Depending on the manner in which he sets up his intuitive concepts of individuals and relations, and expresses his facts.

But what is important to note - for I feel it has caused the greatest confusion - is that sets belong to language, and their corresponding existence in the universe lies in their identifiability , along the lines already indicated, with the entities of the universe.

Such identifiability may be empirically determinable, if the universe is capable of representation as above. Or it may be conventional, when the universe is being constructed.

And, strange as it seems, it is the intuitionist who appears to take the former, Platonistic, stand. And the logicist who possibly adopts the latter.

Instead of trying to name his individuals, relations and facts exhaustively, the logicist only specifies a small central core, and the general rules to which his universe must conform.

Such as his desire that it should contain no indiscernibles. And that every set of individuals is identifiable with a property. And (though, I suspect, under a different guise) that every set of properties defines an individual having just those properties and no others.

And it is on the basis of such desires that the logicist could possibly claim that his language is devoid of content. For if every possible individual and every possible property pertaining to a particular field - characterised by the central core - is to occur in his universe, then, obviously, nothing may be said about what does occur in any particular case. (But this interpretation of no content is likely a far cry from the logicist's actual claim - though the conclusions may be the same.)

Now, conceivably, these rules which the logicist employs for delineating his universe have been misapplied - or they have inherent limitations - as the contradictions indicate.

So it becomes necessary to examine them, and to curtail them. As also to find some possible guidelines for their use and validity. Which is an invitation to the formalist to step in with his apparatus.

And, maybe, every question regarding existence arising from the logicist's activity is - as Carnap suggests - a question of convention as to the type of universe to be delineated.

What position does the so-called Platonic attitude of the logicist then have? And has he not drifted closer to the conventionalist than he would have us believe?

And seeing as how each of the nominalist, intuitionist, logicist, formalist and conventionalist seems to be dealing with a different aspect of the problem, does it make any significance sense to say that they contradict each other?

And could it be that, like the seemingly contradictory wave and particle theories of classical physics, these differing views each contain an important core of truth, hidden by unimportant frills of dispute arising out of possible psychological misunderstandings?

In which case, can I then look forward to a synthesis by extraction - along the lines sketched - of a modern 'wave-mechanical' analogue for mathematical foundations.

### 6. Logicism: Is the Russell-Frege definition of number significant? (29/03/1967)

I do not believe that I have ever been seriously exposed to the influence of nominalism - traditional or Goodman's variety. But I cannot countenance a predicate of predicates unreservedly.

I am able to cheerfully admit the existence of individuals in a universe. I also, hesitantly at first, can embrace the seemingly necessary existence of properties.

But now I see two things.

That each property has an extension, in my language at least, of all the individuals satisfying it. And each individual has an extension of all the properties that it possesses.

And any class of individuals that I am able to construct in my language can only - if at all - be identifiable as the extension of a possible (the job of a formalist being, I believe, to investigate these possibilities - especially if the universe is being set up by convention) property satisfied by the members of the class. The existence of such a property - and hence the reflection of the fact of this existence, in my language - must remain an empirical truth - or a truth by convention.

And, similarly, any class of properties that I can produce in my language is not the reflection of some creature known as a predicate of predicates, but - at the most - the extension identifiable with a possible individual having only the properties contained in the class. The existence of such an individual is again, I dare say, an empirical fact - or a convention.

Now, why does my mind rebel at the thought of indiscriminately creating such individuals?

The reason is chiefly heuristic. As may be expected.

Given a set of individuals, and a two-valued language, I am able to construct $2^n$ distinct classes. If all these exist as properties, then each property is identifiable with some particular class of not more than n individuals. It is not even necessary to insist for the moment that the class be evident to me. So long as I admit that it is a determined class in my language.

Clearly each individual is also identifiable with some class of not more than $2^n$ properties.

But now there are $2^{2^n}$ new individuals which are constructible - at least theoretically so - in my language (which may even embrace a class theory for the construction of its classes - if this is in some way thought possible).

If I try to introduce these in my universe, then the extensions of some of my previous properties will have to be enlarged.

In what sense can I then speak of a property as the static concept it usually is taken to be? Without divorcing it completely from my individuals? In which case, how may I even construct a new property? Unless, of course, I adopt a system of double book-keeping.

And, possibly, this is the reason that Cantor's axiom of comprehension, when applied to ontology, is invalid.

As also the reason that a distinction needs to be drawn between classes and sets in set theory - which is, I believe, implicitly taken to be applicable to both language and ontology. Whether such a distinction has been validly and consistently made relative to the view that I have taken above is a different question. One well worth investigating.

But now I see a major defect in logicism.

$2(f)$ is defined to mean that there exists an $x$, and there exists a $y$, satisfying $f$, and $x$ is not equal to $y$, and if there is some $z$ satisfying $f$, then either $z$ is equal to $x$, or $z$ is equal to $y$.

The class, in my language of course, of $f$'s for which this is true is then identified with an object in the universe containing $f$ over which $x$ and $y$ range.

Such an object, as I have already averred, I can only take to be an individual, say '2'.

But then it appears that every property which has only two true arguments in my universe must necessarily have '2' as one of these (amongst its) arguments! A patently unacceptable conclusion.

At least from an aesthetic point of view, so far as my common sense is concerned. But common sense is not a very reliable guide, and it remains to be seen whether this is also logically (in some sense of the word logic) unacceptable. As I feel it must be. The point is an important one and needs to be investigated.

So I do not accept the individual '2' as identifiable anyhow in my universe. Even though $2(f)$ is a meaningful, and very significant, sentential formula in my language. For it does contain the essence of the meaning-in-use of the number 'two'. And this, I believe, is the really outstanding achievement of logicism. Its analysis of the origin of the number concept. But not its so-called logical construction of the concepts of the integers.

Of course Russell has, to my way of thinking, managed to cloud the issue by ascribing a different level of existence to the individuals constructed from classes of predicates. Which again appears to be a case of multiple standards, since not all classes of predicates - as I have tried to show earlier - need necessarily give rise to the type of difficulty discussed above. Some classes are easily and most naturally identifiable with individuals.

Russell's types are then seen to be nothing more than the setting up of various universes in a kind of chain formation. The lowest being a universe either set up by convention, or which is evident to my senses. The next - not by addition to the first - but rather by identification with expressions of the language in which I talk of my initial universe. And so on.

And of course the language I use to reflect my initial universe will contain expressions for all the possible entities and facts that could possibly occur in it, irrespective of what actually may be occuring at the time I discover/construct it. So Russell may quite readily, though unpardonably for having obfuscated the issue, claim that his universe - which actually contains all the members of the chain that I referred to above - has no content.

And whether we call it one universe or a chain of universes is hardly worth a demonstration at Trafalgar Square, even for Russell.

So long as we can remember that all the successor universes have been constructed from language.

Which gives me enough reason to try and explain why language and ontology have so often been confused.

And my way of justifying the seeming prolificity of language - which I already hinted at above - is this.

I think it would be readily agreed that in the external world there are facts - which may be said to have existence. To ascribe an existence to a non-fact in this universe seems to me somewhat far-fetched, despite McX and Wyman.

Yet I am able, in my language about the external world, to create both factual and non-factual or false expressions.

And this seems a very fortuitious occurence in view of my desire to communicate with, and be communicated to faithfully by, a fallible humanity.

So the expressions in my language seem - at least to my naively finite senses - to exceed the facts in the universe.

Which of course may be an assumption of a very basic and significant nature underlying all my mathematics - hence giving a possible circularity to Cantor's Theorem that $2^n$ exceeds n for all numbers.

Which again needs careful attention.

## 7. Summary (30/03/1967)

01. Discovery of what there 'is', or construction (by convention - other means if thought feasible) of what I feel should be, I take as the basic idea underlying all my mental activity.

02. Language, as the means by which such discovery, or construction, is expressed or conveyed to you.

03. Logical notions as the instruments used to extend what 'is' in any given case to what is possible or could have been possible - in addition to, or as alternative to - the given case.

04. So logic in effect symmetricises language - originally conceived as a carrier of only what there 'is', or, more precisely, of what I believe there 'is' - into containing 'more' than what actually 'is', in terms of what is possible or conceivable.

05. Which gives me a freedom, on the basis of these conceivable entities, entertained by my language (corresponding to the expressions containing free variables, or sets as they are also called) and taking into account what already is, to construct by some means a 'larger', clearly artificial, universe.

06. Larger in the sense that a suitable construction immediately seems to give me Cantor's Theorem - at least if I include all conceivable entities entities of the first into th second.

07. But my constructions necessarily give me a new universe. Though I may be able to map my initial ontology into it in some way.

08. And the obviously recursive procdure gives me a series of universes which Russell calls types.

09. Though there seems no meaningful way in which we can talk of all the universes being united into a universe of universes, with their various entities co-existing peaceably.

10. And the Continuum Hypothesis may be but a convention - a relation between two successive universes - reflecting the manner in which one is constructed out of the other. A relation, then, (like Cantor's) between what is taken 'to be' in a universe, and all that can be constructed from it by means of language.

11. And, so, in some sense what there 'is' does depend on language. At least in all the universes succeeding the initial. And on convention.

12. And the intuitionists have never really gone beyond the first - which rather naturally falls into their domain.

13. While the logicists have seemingly failed to fully appreciate the separate existence of, and distinction between, their various levels. Inspite of the necessity of a special axiom of infinity required at the lowest level.

14. But the nominalist has correctly given the status to what there can be in any universe.

15. And the formalist has justifiably stressed on knowing what probably is, and what could have been in any particular universe. An important question considering our normally incomplete information/construction of any significant initial universe.

16. So all are actually different views of the same thing.

17. And whether this thing is mathematics depends on whether my initial universe has only mathematical entities.

*(Transcribed from original notes: Thursday 8th December 2005 2:17:06 AM by re@alixcomsi.com. Revised: Tuesday 13th December 2005 2:33:26 AM by re@alixcomsi.com.)*